\documentclass[letterpaper, 10 pt, conference]{ieeeconf}  % Comment this line out if you need a4paper

\IEEEoverridecommandlockouts                              % This command is only needed if 
\usepackage{amsmath} % assumes amsmath package installed
\usepackage{tikz}
\usepackage{physics}
\usepackage{subcaption}
\usetikzlibrary{patterns}

\usepackage[version=4]{mhchem}
\usepackage{nomencl}
\usepackage{etoolbox}
\renewcommand\nomgroup[1]{%
  \item[\bfseries
  \ifstrequal{#1}{N}{Notation}{%
  }%
  \ifstrequal{#1}{S}{Sub- \&Superscripts}{%
  }%
  \ifstrequal{#1}{G}{Geometry}{%
  }%
  \ifstrequal{#1}{M}{Math}{%
  }%
  \ifstrequal{#1}{P}{Physics}{%
  }%
]}

\makenomenclature

\title{\LARGE \bf
A Dynamical Simulation Model  of a Cement Clinker Rotary Kiln
}

\author{Jan Lorenz Svensen$^{1,2}$, Wilson Ricardo Leal da Silva$^{2}$, Javier Pigazo Merino$^{2}$,\\ Dinesh Sampath$^{2}$ and John Bagterp J\o rgensen$^{1}$% <-this % stops a space
\thanks{$^{1}$ DTU Compute, Department of Applied Mathematics and Computer Science, Technical University of Denmark, 2800 Lyngby, Denmark
        {\tt\small jlsv@dtu.dk, jbjo@dtu.dk}}%
\thanks{$^{2}$ FLSmdth A/S, 2500, Valby, Denmark
        {\tt\small jls@flsmidth.com, wld@flsmidth.com}}%
}

\begin{document}

\maketitle
\thispagestyle{empty}
\pagestyle{empty}

%%%%%%%%%%%%%%%%%%%%%%%%%%%%%%%%%%%%%%%%%%%%%%%%%%%%%%%%%%%%%%%%%%%%%%%%%%%%%%%%
\begin{abstract}
 This study provides a systematic description and results of a dynamical simulation model of a rotary kiln for clinker, based on first engineering principles.  The model is built upon thermophysical, chemical, and transportation models for both the formation of clinker phases and fuel combustion in the kiln. The model is presented as a 1D model with counter-flow between gas and clinker phases and is demonstrated by a simulation using industrially relevant input. An advantage of the proposed model is that it provides the evolution of the individual compounds for both the fuel and clinker. As such, the model comprises a stepping stone for evaluating the development of process control systems for existing cement plants.

\end{abstract}
%%%%%%%%%%%%%%%%%%%%%%%%%%%%%%%%%%%%%%%%%%%%%%%%%%%%%%%%%%%%%%%%%%%%%%%%%%%%%%%%
\section{Introduction}
In this paper, a differential-algebraic model of the dynamics of a rotary kiln is proposed, with the intention of modeling the formation of clinker and combustion.
The modeling approach is primarily based on first engineering principles, to provide an intuitive understanding of the processes included.% within the model.

\subsection{Background and motivation}
The cement industry accounts for 8\% of the world's \ce{CO_2} emissions \cite{CO2Techreport}. A main contributor to these \ce{CO_2} emissions and energy consumption is the production of clinker; a main component in the cement.
For example, Ordinary Portland Cement comprises ca. 90\% clinker, with the remainder being gypsum and other mineral additions (e.g. limestone).
%e.g. Ordinary Portland Cement is more than $90\%$ clinker (the clinker ratio), with the rest being gypsum.
The strength of cement depends on its composition, especially the amount of \ce{C_3S}, and thus depends on the composition of the clinker, the clinker quality, to provide enough strength.
%If the clinker quality is uncertain or has variation, then to ensure the suitable cement strength a higher clinker ratio is used to guarantee enough of each compound.
%To improve the sustainability of cement, reducing uncertainty and variation of the clinker quality is thus needed, to reduce the need for clinker and thus \ce{CO2} emissions. 
In addition, considerable fluctuations in the clinker composition can limit the utilization of supplementary cementitious materials (SCMs) - because a greater percentage of clinker is likely required to cover for the resulting variability in cement strength.
In view of that, reducing uncertainty and variations of the clinker quality is thus a must to enable a more sustainable and controlled cement production.

The uncertainty can be addressed by improving the quality predictions and the variation through control of the operation. The design of both approaches requires the availability of system responses, from either an actual system or a simulation model, to account for the dynamical transitions between operation states. 
Although a simulation model does not account for all process details, it enables cost-effective and faster testing without disrupting a production site.

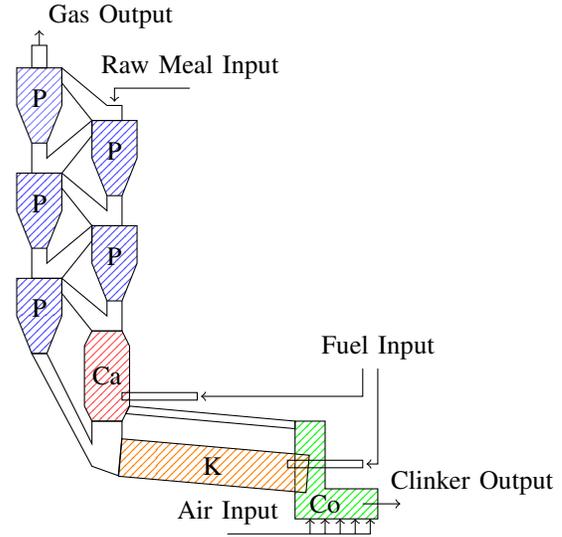
\begin{figure}[tb]
    \centering
        \begin{tikzpicture}    
            \draw[pattern=north east lines, pattern color=blue!50] (0.2,5.5) -- (0.2,6) -- (0.8,6) -- (0.8,5.5) -- (0.6,5) -- (0.4,5)-- cycle;
            \draw[pattern=north east lines, pattern color=blue!70] (1.2,4.8) -- (1.2,5.3) -- (1.8,5.3) -- (1.8,4.8) -- (1.6,4.3) -- (1.4,4.3)-- cycle;
            \draw[pattern=north east lines, pattern color=blue!70] (0.2,4.1) -- (0.2,4.6) -- (0.8,4.6) -- (0.8,4.1) -- (0.6,3.6) -- (0.4,3.6)-- cycle;
            \draw[pattern=north east lines, pattern color=blue!70] (1.2,3.4) -- (1.2,3.9) -- (1.8,3.9) -- (1.8,3.4) -- (1.6,2.9) -- (1.4,2.9)-- cycle;
            \draw[pattern=north east lines, pattern color=blue!70] (0.2,2.7) -- (0.2,3.2) -- (0.8,3.2) -- (0.8,2.7) -- (0.6,2.2) -- (0.4,2.2)-- cycle;

            \draw  (1.4,5.5) -- (0.8,6) -- (0.8,5.8) -- (1.2,5.3)-- (1.6,5.3) -- (1.6,5.5) -- cycle;
            \draw (0.6,5) -- (0.6,4.8) -- (1.2,5.3) -- (1.2,5.1) -- (0.8,4.6)-- (0.4,4.6) -- (0.4,5);
            \draw (1.4,4.3) -- (1.4,4.1) -- (0.8,4.6) -- (0.8,4.4) -- (1.2,3.9)-- (1.6,3.9) -- (1.6,4.3);
            \draw (0.6,3.6) -- (0.6,3.4) -- (1.2,3.9) -- (1.2,3.7) -- (0.8,3.2)-- (0.4,3.2) -- (0.4,3.6);
            \draw (1.4,2.9) -- (1.4,2.7) -- (0.8,3.2) -- (0.8,3.0) -- (1.2,2.5)-- (1.6,2.5) -- (1.6,2.9);
            \draw (0.5,2.8) node{P};
            \draw (1.5,3.5) node{P};
            \draw (0.5,4.2) node{P};
            \draw (1.5,4.9) node{P};
            \draw (0.5,5.6) node{P};
            
            % Calciner
            \draw[pattern=north east lines, pattern color=red!70] (1.1,1.5) -- (1.1,2.3) -- (1.2,2.5) -- (1.6,2.5) -- (1.7,2.3) -- (1.7,1.5) -- (1.6,1.3) -- (1.2,1.3)-- cycle;
            \draw (1.4,1.9) node{Ca};
            
            % kiln and cooler
            \draw[pattern=north east lines, pattern color=green] (3.9,1+0.3) -- (3.9,1-1) -- (5,1-1) -- (5,1-0.6) -- (4.3,1-0.6) -- (4.3,1+0.3)-- cycle;
            \draw[rotate around={-5:(0,0)},pattern=north east lines, pattern color=orange] (1.5,1-0.3) -- (4,1-0.3) -- (4,1+0.2) -- (1.5,1+0.2) -- cycle;
            \draw (2.8,0.7) node{K};
            
            \draw (3.9,1+0.3) -- (3.9,1+0.2) -- (1.65,1.4) -- (1.7,1.5) -- cycle;
            \draw (1.6,1.3) -- (1.6,1.1) -- (1.55,0.575) -- (1.2,0.7) -- (0.4,2.2) -- (0.6,2.2) -- (1.2, 1) -- (1.2,1.3) -- cycle;
            \draw (4.3,0.2) node{Co};

            % fuel and gas pipes in/out
            \draw  (3.8,0.68) rectangle (4.8,0.78);
            \draw  (1.6,1.58) rectangle (2.6,1.68);
            \draw  (0.4,6) rectangle (0.6,6.3);

            \draw[->] (4.8,0.2) -- node[anchor=south west]{Clinker Output} (5.3,0.2);
            \draw[->] (5.0,2) node[anchor=south]{Fuel Input} -- (5.0,0.73)  -- (4.85,0.73);
            \draw[->] (4.8,2) -- (4.8,1.63)  -- (2.65,1.63);
            \draw[->] (0.5,6.3) -- node[anchor=south west]{Gas Output} (0.5,6.5);
            \draw[->] (2.5,5.73) node[anchor=south]{Raw Meal Input} -- (1.5,5.73)  -- (1.5,5.55);
            \draw[->] (3.0,-0.2) node[anchor=south]{Air Input} -- (4.1,-0.2)  -- (4.1,0.0);
            \draw[->] (4.1,-0.2) -- (4.3,-0.2)  -- (4.3,0.0);
            \draw[->] (4.3,-0.2) -- (4.5,-0.2)  -- (4.5,0.0);
            \draw[->] (4.5,-0.2) -- (4.7,-0.2)  -- (4.7,0.0);
            \draw[->] (4.7,-0.2) -- (4.9,-0.2)  -- (4.9,0.0);
    \end{tikzpicture}    
    \caption{Simplified Clinker Production Line: preheating Cyclones ({\color{blue}P}), Calciner ({\color{red}Ca}), Rotary Kiln ({\color{orange}K}), and Cooler ({\color{green}Co})}\label{fig:production}
\end{figure}

The production of clinker comprises several process equipment: a pre-heating tower (with a set of cyclones), a calciner, a rotary kiln, and a cooler, see Fig. \ref{fig:production}. 
As the formation of clinker takes place in the rotary kiln, this paper focuses on formulating a dynamical model for the rotary kiln.

\subsection{Research review}
Several models have been proposed in the literature to describe the kiln. Most of these focuses on steady-state descriptions of the kiln \cite{Mujumdar2006}-\nocite{Hanein2017}\cite{Mastorakos1999CFDPF}, while few give a dynamic description \cite{Spang}-\nocite{SUN2020,LIU2016}\cite{GINSBERG2011}.

Spang's dynamic model \cite{Spang} includes a description of the bed chemistry and temperature energy balances. In Sun's model \cite{SUN2020}, the temperature is also described, while the mass balance only considers the overall bulk matter. Similarly, in Liu's model, \cite{LIU2016} the bulk bed/gas is considered with added momentum balance; while Ginsberg \cite{GINSBERG2011} assumes a steady-state gas balance in its dynamic description.

Mastorakos \cite{Mastorakos1999CFDPF} presented a detailed steady-state model with a 3D temperature energy balance, a 1D mass balance, and constant velocities. The combustion and gas parts were only mentioned briefly in words.
Further, the steady-state model of Hanein \cite{Hanein2017} only contains the energy balance, while Mujumdar's model \cite{Mujumdar2006} outlines the temperature energy balances and solid mass balance accounting for particle size, though it does not consider the gas mass balance nor the fuel part.

% \cite{ARIYARATNE2015} particle fuel dynamics
%This paper aims to present a dynamic model, which can be described in a simple way: that preserves the understanding of each part of the physical system, and allow for easy extension/modification to fit new scenarios.

\subsection{Contribution}
In none of the reported models, the fuel chemistry is accounted for, and only few consider the gas mass dynamics and velocity changes. Thus, it is not trivial to potentially include modern considerations, e.g. the effects of ash content from alternative fuels or gas reactions giving non-\ce{CO2} emissions e.g. \ce{NO_x}.

To address these limitations, this paper aims to present a dynamic kiln model where: 1) the description preserves the physical overview and understanding of the system; 2) it is easy to modify, for the inclusion of new scenarios or descriptions; and 3) the dynamics of the individual compounds is modeled, e.g. \ce{CO2} content.

The presented model covers the main dynamics of the kiln with the following simplifications and assumptions: 1) it is a 1D-model (axial); 2) the heat loss to the environment is negligible; and 3) only the primary clinker- and combustion reactions are included.

The standard cement chemist notation will be used for the following compounds, i.e.:
 \ce{(CaO)_2SiO_2} as \ce{C_2S}, \ce{(CaO)_3SiO_2} as \ce{C_3S}, \ce{(CaO)_3Al_2O_3} as \ce{C_3A} and \ce{(CaO)_4(Al_2O_3)(Fe_2O_3)} as \ce{C_4AF}, where C, A, S and F refers to \ce{CaO}, \ce{Al_2O_3}, \ce{SiO_2} and \ce{Fe_2O_3}, respectively. Finally $\partial_x$ will note the differential $\pdv{}{x}$

\section{Model Layout}
The kiln model is formulated as a Differential-Algebraic system, with
dynamic descriptions of the concentrations $C$ and internal energy densities $\hat U$, and algebraic relations for the temperatures $T$, pressures $P$, and fill angle $\theta$:
\begin{subequations}\label{eq:dyn}
\begin{align}
    \partial_tx& = f(x,y),\quad x=[C;\hat U]\\
    0 &= g(x,y),\quad y=[T;P;\theta].
\end{align}
\end{subequations}
%The solution for \eqref{eq:dyn} can be computed numerically using e.g. implicit Euler's method and Newton-Raphson's method. 

For clarity, the different aspects of the kiln model are discussed separately in order of computation: A) Thermophysical model, B) Geometry, C) Transportation, D) Stoichiometry \& Kinetics, E) Mass balances, F) Energy balances, and G) Algebraic relations.

%To ease the computation, clarity and understanding, the different aspects of the kiln model are kept separate. The model consists of the following parts in order of computation: A) Thermodynamics, B) Geometry, C) Transportation, D) Stoichiometry \& Kinetics, E) Mass balances, F) Energy balances, and G) Algebraic relations.

\section{Kiln Model}
A rotary kiln is a rotating cylinder of length $L$ and inner radius $r_c$ with an inclination $\psi$, Fig. \ref{fig:profiles} illustrates the kiln's cross and axial profiles. 
We use a finite-volume approach to describe the kiln in $n_v$ segments of length $\Delta z = L/n_v$ and volumes $V_{\Delta} = \pi r_c^2\Delta z$.
The model segment consists of 3 phases; wall ($_w$), solids ($_s$), and gasses ($_g$). The concentrations are of the individual compounds of each phase.%, to allow for reaction kinetics.

% {\color{red}The states in the kiln model are the concentration profiles $C$ for each substance and the energy density profiles $\hat{U}$ of the wall ($_w$), solid ($_s$), and gas ($_g$), along the longitudinal axis of the kiln. 
% The known initial conditions include Similar profiles of the temperatures $T$ and pressures $P$.}

\subsection{Thermophysical model}
If we define concentration $C$ as mole per total segment volume $V_{\Delta}$ and assume gasses are ideal, then the enthalpy and volumes for the solid and gas part can be defined as functions $H$ and $V$, linear w.r.t. the moles, e.g.:
\begin{align}
    H(T,P,n) &= \sum_i (\Delta H_{\mathrm{f},i}^\circ + \int^T_{T_0} c_{m,i}(\tau) d\tau ) n_i\\
    V_s(T,P,n) &= \sum_i \frac{M_i}{\rho_i}n_i,\quad  V_g(T,P,n) = \frac{RT}{P}\sum_i n_i.
\end{align}

The enthalpy- and volume densities are obtained by the H and V functions, by applying the concentrations instead:
%If we apply the concentrations instead, the H and V functions can define the enthalpy- and volume densities as:
\begin{align}
    \hat{H}_s &= H(T_s,P,C_s),\quad  &\hat{V}_s &= V(T_s,P,C_s)\\ %&= \sum_i H_i(T,P,n_i) =\sum_i h_i(T,P) n_i \\
    \hat{H}_g &= H(T_g,P,C_g),\quad &\hat{V}_g &= V(T_g,P,C_g). %&= \sum_i V_i(T,P,n_i) = \sum_i \frac{M_i}{\rho_i} n_i
\end{align}
with the solid- and gas volumes obtained by their densities:
\begin{align}
    V_g = \hat{V}_gV_{\Delta},\quad  V_s = \hat{V}_sV_{\Delta}.
\end{align}
The energy densities of each phase are then described by:
\begin{align}
    \hat{U}_g = \hat{H}_g - P\hat{V}_g, \quad \hat{U}_s = \hat{H}_s, \quad \hat{U}_w = H(T_w,P,\rho_w).\label{eq:EnergyAlgebra}
\end{align}

% We define concentration $C$ as mole per total segment volume $V_{\Delta}$ and assume gasses are ideal. For each segment, the enthalpy and volumes of the solid and gas part read:
% \begin{align}
%     H_s &= H(T_s,P,n_s),\quad  & V_s &= V(T_s,P,n_s)\\ %&= \sum_i H_i(T,P,n_i) =\sum_i h_i(T,P) n_i \\
%     H_g &= H(T_g,P,n_g),\quad &V_g &= V(T_g,P,n_g) %&= \sum_i V_i(T,P,n_i) = \sum_i \frac{M_i}{\rho_i} n_i 
% \end{align}
% with the functions $H$ and $V$ being linear w.r.t. the moles.
% $H(T,P,n) = \sum_i\int cp_i(\tau) d\tau|_{\tau=T} n_i$

% $V(T,P,n) = \sum_i M_i/\rho_i n_i , RT/P\sum_i n_i$

% If we apply the concentrations instead, the H and V functions can define the enthalpy- and volume densities as:
% \begin{align}
%     \hat{H}_s &= H(T_s,P,C_s),\quad  &\hat{V}_s &= V(T_s,P,C_s)\\ %&= \sum_i H_i(T,P,n_i) =\sum_i h_i(T,P) n_i \\
%     \hat{H}_g &= H(T_g,P,C_g),\quad &\hat{V}_g &= V(T_g,P,C_g). %&= \sum_i V_i(T,P,n_i) = \sum_i \frac{M_i}{\rho_i} n_i .
% \end{align}
% The internal energy densities of the gas, solid, and wall parts can be described by:
% \begin{align}
%     \hat{U}_g = \hat{H}_g - P\hat{V}_g, \quad \hat{U}_s = \hat{H}_s, \quad \hat{U}_w = \hat{H}_w(T_w).\label{eq:EnergyAlgebra}
% \end{align}
% The solid and gas volumes can be obtained from their densities:
% \begin{align}
%     V_g = \hat{V}_gV_{\Delta},\quad  V_s = \hat{V}_sV_{\Delta}.
% \end{align}

\subsection{Geometric relationships}

\begin{figure}
    \centering
    \begin{subfigure}[b]{0.5\textwidth}
    \centering
        \begin{tikzpicture}        
            \draw[pattern=north east lines, pattern color=brown!50] (2,2) circle (2 cm);
            \draw[thick,-] (2,2) circle (1.5 cm);
            \fill[white] (2,2) circle (1.5 cm);
            \fill[blue!40!white] (2,2) circle (0.1 cm);
            \draw[pattern=north west lines, pattern color=orange!60]  (1.5,0.6) arc (-110:0:1.5cm)  -- cycle;
    
            \draw[ultra thick, ->] (2,4.2)  arc (90:150:2.2cm) node[anchor=south east]{$\omega$};\draw[] (2,4.4);
            \draw[thick,<->] (2,2) -- node[anchor=south]{$r_c$} (3.05,3.05) ;
    
            \draw[thick,<->] (1.5,0.6) -- node[anchor=south]{$L_c$} (3.5,2) ;        
            \draw[-] (1.5,0.6) -- (2,2) ;
            \draw[-] (3.5,2) --  (2,2) ;
            \draw (2.2,2)node[anchor=north]{$\theta$} arc (0:-112:0.2cm) ;
            %\draw[thick,<->] (2.5,1.3) -- node[anchor=south]{$L_c$} (2,2) ;
            \draw[thick,<->] (2.5,1.3) node[anchor=west]{$h$} --  (2.88,0.775) ;
    
            \draw[dashed,-] (1.5,0.6) -- (0.65,0) ;
            \draw[dashed,-] (2,0) -- (0.65,0) ;
            \draw (1.0,0)node[anchor=north west]{$\xi$} arc (0:35:0.35cm) ;
    
            \draw[-] (3.7,2.5) -- (4,3.0) node[anchor=south]{$A_w$};
            \draw[-] (3.2,1.5) -- (4.5,1.5) node[anchor=south]{$A_s$};
            \draw[-] (2.0,3.0) -- (3.5,4.0) node[anchor=west]{$A_g$};
    
            \draw[thick,->] (1.1,2.5) -- node[anchor=west]{$Q_{gw}$} (0.5,3.0) ;
            \draw[thick,->] (1.8,1.2) -- node[anchor=east]{$Q_{gs}$} (2.1,0.8) ;
            \draw[thick,->] (2,0.3) -- node[anchor=west]{$Q_{ws}$} (2,0.7) ;
        \end{tikzpicture}    
        \caption{Cross section.}
        \label{fig:cross}
    \end{subfigure} 
    
    \begin{subfigure}[b]{0.5\textwidth}
    \centering
    \resizebox{1\textwidth}{!}{%
    \begin{tikzpicture}   
        \draw[thick,->] (0,-0.6)node[anchor=north]{$0$} -- (8.7,-0.6)node[anchor=north]{$L$}node[anchor=west]{$z$};
        \draw[pattern=north east lines, pattern color=brown!50] (0,0.0) rectangle (8.7,4.0);
        \fill[white] (0,0.5) rectangle (8.7,3.5);
        \draw (0,0.5) rectangle (8.7,3.5);

        \draw[pattern=north west lines, pattern color=orange!60]  (0,1.9) -- (8.7,1.2) -- (8.7,0.5) -- (0,0.5)-- cycle;
        \draw (8.95,0.0625) -- (6.7,-0.5); 
        \draw  (7.3,-0.3) arc (-112:-180:0.3cm) node[anchor= north east]{$\psi$};

        \draw[dashed,blue!40!white] (0,2) -- (9,2);
        \draw[dashed,blue!40!white] (6.9,1.4) -- (8.3,1.4);
        \draw (8.2,1.4)node[anchor=west]{$\phi$} arc (0:-35:0.25cm) ;

        \draw (2.5,-0) rectangle node[anchor=south]{$z_{k}$}(4.0,4);
        \draw (2.5,-0.7) -- (2.5,-0.5) node[anchor=south]{$z_{k-\frac{1}{2}}$};
        \draw (4.0,0) rectangle node[anchor=south]{$z_{k+1}$} (5.5,4) ;
        \draw (4.0,-0.7) -- (4.0,-0.5) node[anchor=south]{$z_{k+\frac{1}{2}}$};
        \draw (5.5,-0.7) -- (5.5,-0.5) node[anchor=south]{$z_{k+\frac{3}{2}}$};
        \draw[<->] (2.5,-0.7) -- node[anchor=north]{$\Delta z$} (4.0,-0.7);
        \draw[<->] (4.0,-0.7) -- node[anchor=north]{$\Delta z$} (5.5,-0.7);
        
        %\draw[thick,<->] (2.5,0.5)  --node[anchor=east]{$h$}  (2.5,1.5) ;
        \draw[thick,<->] (4.0,0.5)  --node[anchor=west]{$h_{k+\frac{1}{2}}$}  (4.0,1.4) ;
        %\draw[thick,<->] (5.5,0.5)  --node[anchor=east]{$h$}  (5.5,1.25) ;

        \draw[thick,<-] (2.25,2.8) -- node[anchor=north east]{$\Tilde{H}_{g,k-\frac{1}{2}}$} node[anchor=south east]{$N_{g,k-\frac{1}{2}}$} (2.75,2.8);
        \draw[thick,<-] (5.25,2.8) -- node[anchor=south west]{$N_{g,k+\frac{3}{2}}$}node[anchor=north west]{$\Tilde{H}_{g,k+\frac{3}{2}}$} (5.75,2.8);
        \draw[thick,->] (2.25,1.1)node[anchor=east]{\small$\Tilde{H}_{s,k-\frac{1}{2}}$} node[anchor=north east]{\small$N_{s,k-\frac{1}{2}}$}node[anchor=south east]{\small$\Tilde{Q}_{s,k-\frac{1}{2}}$} --  (2.75,1.1);
        \draw[thick,->] (5.25,1.0) -- node[anchor=north west]{\small$N_{s,k+\frac{3}{2}}$} node[anchor=south west]{\small$\Tilde{Q}_{s,k+\frac{3}{2}}$} (5.75,1.0)node[anchor=west]{\small$\Tilde{H}_{s,k+\frac{3}{2}}$};
        \draw[thick,<->] (3.2,1.4) --  (3.2,1.8)node[anchor= west]{$R_{k}$};

        \node[anchor=south west,inner sep=0] at (0.1,0.1){Wall};
        \node[anchor=south west,inner sep=0] at (0.1,3.6){Wall};
        \node[anchor=south west,inner sep=0] at (0.1,1.1){Solid};
        \node[anchor=south west,inner sep=0] at (0.1,3.1){Gas};
    \end{tikzpicture}
    }
    \caption{Axial profile with segment notation for the $k$-th volume.}
    \label{fig:axial}
    \end{subfigure}
    \caption{Kiln profiles}\label{fig:profiles}
\end{figure}

%%%%%  General        %%%%%
The gas and solid volumes across the kiln define the geometric aspect of each segment. The solid volume $V_s$ gives the bed height $h$, through the cross-section area:
\begin{align}
    V_s(z_k) = \int_{z_{k-\frac{1}{2}}}^{z_{k+\frac{1}{2}}} A_s(z) dz.
\end{align}
As the cross section $A_s$ of the solids forms a circle segment \cite{CircSeg}, see Fig. \ref{fig:cross}, the cross-area defines the
bed height $h$, the surface chord $L_c$, the fill angle $\theta$, and the bed slope angle $\phi$: 
% \begin{align}
%     A_s(h) &= r_c^2 acos\bigg(\frac{r_c-h}{r_c}\bigg) - (r_c-h)\sqrt{2r_ch-h^2} \label{eq:geoAlg}\\
%     L_c(h) &= 2 \sqrt{2r_ch-h^2}\\
%     \theta(h) &= 2 asin\bigg(\frac{L_c(h)}{2r_c}\bigg),\quad 
%     \phi(h) = atan\bigg(-\pdv{h}{z}\bigg).
% \end{align}
\begin{align}
    A_s(z) &= \frac{r_c^2}{2}(\theta(z) - \text{sin}(\theta(z)) )\label{eq:geoAlg}\\
    L_c(z) &= 2 r_c\text{sin}(\frac{\theta(z)}{2}),\quad h(z) = r_c(1-\text{cos}(\frac{\theta(z)}{2}))\\          
    \phi(z) &= \text{atan}\bigg(-\pdv{h(z)}{z}\bigg).
\end{align}
The surface areas between the wall, gas, and solid along the segment, and the gas cross area $A_g$ are given by: 
\begin{align}
    A_g(z) & = \pi r_c^2 - A_s(z)\\
    A_{gs}(z_k) &= \int_{z_{k-\frac{1}{2}}}^{z_{k+\frac{1}{2}}} L_c(z) dz\\
    A_{ws}(z_k) &= \int_{z_{k-\frac{1}{2}}}^{z_{k+\frac{1}{2}}} r_c\theta(z) dz\\
    A_{gw}(z_k) &= 2\pi r_c\Delta z - A_{ws}(z_k)
\end{align}
where $A_{gs}$ is the gas-solid surface area, $A_{gw}$ is the gas-wall surface area, and $A_{ws}$ is the wall-solid surface area.

\subsection{Transportation}
The transportation of matter and energy in the kiln depends on material flux, heat flux, heat convection, and heat radiation. The details on these aspects are described as follows.

\subsubsection{Material flux} The flux of each compound is defined independently for solid and gas  given their opposite flow. It consists of an advection term and a diffusion term:
\begin{align}
    N_{i,s} &= v_{s}C_{i,s} + J_{i,s},\quad &J_{i,s}& = -D_{i,s}\partial_zC_{i,s}\\
    N_{i,g} &= v_{g}C_{i,g} + J_{i,g},\quad &J_{i,g}& = -D_{i,g}\partial_zC_{i,g}
\end{align}
with the diffusion $J_{i}$ given by Fick's law \cite{Perry}.
\paragraph{Velocities}
The advection term depends on the concentration and the phases' velocity. The velocity of the solid depends on the rotational velocity of the kiln $\omega$:
\begin{align}
    v_{s}(z) = \omega \frac{\psi+\phi(z) \text{cos}(\xi)}{\text{sin}(\xi)}\frac{\pi L_c(z) }{\text{asin}(\frac{L_c(z)}{2r_c})}
\end{align}
$v_s$ is the average velocity at the location, %defined by the formula for the average axial transport velocity of a particle in the bed 
accounting for the cascading motion due to the rotation \cite{Saeman1951}. The repose angle $\xi$ is related to the rotational velocity by
\begin{equation}
    \xi = a_\omega\omega + b_\omega
\end{equation}    
with the coefficients determined experimentally \cite{Yang2003,Yamane1998}.
% \begin{align}
%     v_{s} = 2\pi\omega\frac{\psi+\phi cos(\xi)}{sin(\xi)}\frac{2r_c-h_{z} }{2}
% \end{align}

Given that the gas velocity in a kiln is below 0.2 Mach speed \cite{MortenPHD}, according to Howel \& Weathers \cite{Darcy-Howel}, the Darcy-Weisbach equation in (\ref{DWeq}) is valid to describe the gas velocity, despite the gas being compressible.
\begin{align} 
    v_g^2 = 2\frac{|\Delta P|}{\Delta z}\frac{D_H}{f_D\rho_g}, \label{DWeq}
\end{align}
Assuming turbulent flow, the Darcy friction factor $f_D$ is% given by
\begin{align}
    f_D = \frac{0.316}{\sqrt[4]{Re}} = \frac{0.316\sqrt[4]{\mu_g}}{\sqrt[4]{\rho_g |v_g|D_H}}, \ D_H = \frac{4V_g}{A_{gw}+A_{gs}}
\end{align}
%laminar flow, $f_D$ is given by
% \begin{align}
%     f_D = \frac{64}{Re} = \frac{64\mu_g}{\rho_g v_g D_H}, \quad D_H = \frac{V_g}{A_{gw}+A_{gs}}
%\end{align}
where $D_H$ is the hydraulic diameter for a Non-uniform nor circular channel \cite{HESSELGREAVES20171}.
The gas velocity is then given by
\begin{align} 
    v_g = \Big(\frac{2}{0.316}\sqrt[4]{\frac{D_H^{5}}{\mu_g\rho_g^3}}\frac{|\Delta P|}{\Delta z}\Big)^{\frac{4}{7}}\text{sgn}\Big(\frac{\Delta P}{\Delta z}\Big),\label{eq:vel}
\end{align}
% \begin{align} 
%     v_g = \frac{D_H^2}{32\mu_g }\frac{\Delta P}{\Delta z}
%\end{align}
since negative pressure difference leads to a negative flow.

\paragraph{Diffusion}
The diffusion of each phase is defined by the diffusion coefficients $D_{i,s}$ and $D_{i,g}$. 
% {\color{red}For the solid case, the coefficient is assumed to be governed by Arrhenius equation:
% \begin{equation}
%      D_{i,s} = D_0e^{-\frac{E_{A}}{RT}}
% \end{equation}
% where the $D_0$ coefficient  and activation energy $E_A$ are determined experimentally.}
According to Mujumdar \cite{Mujumdar2006}, solid diffusion can be neglected since an industrial kiln has a Peclet number greater than $10^4$, $D_s=0$.

For diffusion within a gas mixture, the coefficient is given by (\ref{eq:gdiff}-\ref{eq:gdiff2}) using Fuller's model \cite{Poling2001Book}, where $x_j$ is the mole fraction, $V_{\sum}$ the diffusion volume, and $M$ the molar mass.
\begin{align}    
     D_{i,g} & = \bigg(\sum_{\substack{j=1\\j\neq i}}\frac{x_j}{c_gD_{ij}}\bigg)^{-1},\ x_j = \frac{C_{j,g}}{c_g}\label{eq:gdiff}, \ c_g= \sum_jC_{j,g}\\
    D_{i,j} &= \frac{0.00143T^{1.75}}{P M_{ij}^{\frac{1}{2}}[(V_{\sum})^{\frac{1}{3}}_i+(V_{\sum})^{\frac{1}{3}}_j]^2}\label{eq:gdiff2}, \ M_{ij} = \frac{2}{\frac{1}{M_i}+\frac{1}{M_j}} 
\end{align}

\subsubsection{Heat conduction}
The heat flux of each phase is given by Fourier's law \cite{Perry}:
\begin{align}
    \Tilde{Q}_i = -k_i\partial_z{T_i}, \quad \text{for }i=\{s,g,w\}
\end{align}
with $k_i$ being the conductivity of phase \textit{i}, e.g. $k_s$ for the solid mixture.

\subsubsection{Heat convection}
The transfer of heat due to convection between the gas, solid, and wall is given by \cite{Mujumdar2006}:
\begin{align}
    Q_{gs}^{cv} &= A_{gs}\beta_{gs}(T_g-T_s)\label{eq:gasconv1}\\
    Q_{gw}^{cv} &= A_{gw}\beta_{gw}(T_g-T_w)\label{eq:gasconv2}\\
    Q_{ws}^{cv} &= A_{ws}\beta_{gw}(T_w-T_s)
\end{align}
where $A_{ij}$ is the in-between surface area, and $\beta_{ij}$ is the convection coefficient. The coefficients of the three phases are defined by Tscheng \cite{Tscheng1979} as:
\begin{align}
    \beta_{gs} =& \frac{k_g}{D_e}Nu_{gs}, \ Nu_{gs}=0.46Re_D^{0.535}Re_\omega^{0.104}\eta^{-0.341}\\
    \beta_{gw} &= \frac{k_g}{D_e}Nu_{gw}, \quad Nu_{gw}=1.54Re_D^{0.575}Re_\omega^{-0.292}\\
    \beta_{ws} &= 11.6\frac{k_s}{l_w}(\frac{\omega r_c^2\theta}{\alpha_B})^{0.3}, \ \alpha_B = \frac{k_s}{\rho_sC_{ps}},\ l_w = r_c\theta
\end{align}
where $k_i$ is the conductivity of phase i, $D_e$ the effective diameter, $\alpha_B$ the thermal diffusivity, and $l_w$ the contact perimeter between the solid and wall.  $Re_D$ and $Re_\omega$ are Reynold's numbers and $\eta$ is the solid fill fraction:
\begin{align}
Re_D &= \frac{\rho_gv_gD_e}{\mu_g}, \quad Re_\omega = \frac{\rho_g\omega D_e^2}{\mu_g}\\
    \eta &= \frac{\theta-\text{sin}(\theta)}{2\pi}\quad
D_e = 2r_c \frac{\pi-\frac{\theta}{2}+\frac{\text{sin}(\theta)}{2}}{\pi-\frac{\theta}{2}+\text{sin}(\frac{\theta}{2})}
\end{align}
with $\mu_g$ being the viscosity of the gas mixture.

The densities $\rho$ and heat capacity $C_{ps}$ are given by:
\begin{align}
    \rho_j = \frac{\sum_iM_iC_{j,i}}{\hat{V}_j},\quad C_{ps} = \sum_in_ic_{m,i}
\end{align}
with $M$ being the molar mass and $c_m$ the molar heat capacity.

\paragraph{Viscosity \& conductivity} 
For a single gas, Sutherland's formula for temperature-dependent viscosity read \cite{Sutherland1893}:
\begin{align}
    \mu_{g,i} = \mu_0 \bigg(\frac{T}{T_0}\bigg)^{\frac{3}{2}}\frac{T_0+S_\mu}{T+S_\mu},
\end{align}
where $S_\mu$ can be calibrated given two measures of viscosity.

The viscosity and conductivity of a mixed gas are given by the similar formulas of Wilke \cite{Wilke1950} and Wassiljewa \cite{Poling2001Book}:
\begin{align}
    \mu_g &= \sum_i\frac{x_i\mu_{g,i}}{\sum_jx_j\phi_{ij}},\quad k_g = \sum_i\frac{x_ik_{g,i}}{\sum_jx_j\phi_{ij}}\\
    \phi_{ij} &= \bigg(1+\sqrt{\frac{\mu_{g,i}}{\mu_{g,j}}}\sqrt[4]{\frac{M_j}{M_i}}\bigg)^2\bigg(2\sqrt{2}\sqrt{1+\frac{M_i}{M_j}}\bigg)^{-1}.
\end{align}
%with $x_i$ being the mole fraction.

%For the conductivity of mixed solids, no formula was found in the literature, so we apply the formula for conductivity through layers \cite{Perry}, adapted for volumetric ratios:
For the conductivity of mixed solids, we apply the formula for conductivity through layers \cite{Perry} using volumetric ratios, assuming layers with length $\Delta z_i$ and constant cross area:
\begin{align}
    k_s = \frac{1}{\sum_i\frac{V_{s,i}}{V_{s}}\frac{1}{k_{s,i}}}, \quad \frac{V_{s,i}}{V_{s}} = \frac{\Delta z_i}{\Delta z}.
\end{align}

% {\color{red}No formula was found in the literature for the conductivity of mixed solids, so inspired by the formula for thermal conduction through layers \cite{Perry}, where the weight depends on the total distance ratio with a given conductivity, the following formula is utilized:
% \begin{align}
%     k_s = \frac{1}{\sum_i\frac{V_{s,i}}{V_{s}}\frac{1}{k_{s,i}}}
% \end{align}
% where the weight has been reinterpreted as a volume ratio. }

\subsubsection{Heat Radiation}
Assuming axial radiation is negligible, then the heat transfer due to radiation is given by \cite{Mujumdar2006}:
\begin{align}
    Q_{gs}^{rad} &= \sigma A_{gs}(1+\epsilon_s)\frac{\epsilon_gT^4_g-\alpha_gT^4_s}{2}\\
    Q_{gw}^{rad} &= \sigma A_{gw}(1+\epsilon_w)\frac{\epsilon_gT^4_g-\alpha_gT^4_w}{2}\\
    Q_{ws}^{rad} &= \sigma A_{ws}\epsilon_w\epsilon_s\Omega(T^4_w-T^4_s),\quad 
    \Omega = \frac{L_c}{2(\pi-\psi)r_c},
\end{align}
where $\sigma$ is Stefan-Boltzmann's constant, $\epsilon$ is the emissivity, $\alpha$ is the absorptivity, and $\Omega$ is the form factor.

\paragraph{Emissivity \& absorptivity}
In the literature, the standard values of the emissivity for the kiln wall and the solids are $\epsilon_w=0.85$ and $\epsilon_s=0.9$ respectively \cite{Hanein2017}.
The emissivity of the gas $\epsilon_g$ can be computed using the WSGG model of 4 grey gases \cite{Johanson2011} as shown below, assuming only \ce{H2O} and \ce{CO2} affect the radiation, with other gasses being transparent.
\begin{align}
    \epsilon_g &= \sum^4_{j=0}a_j(1-e^{-k_jS_mP(x_{H2O}+x_{CO2})})\\
    a_0 &= 1-\sum^4_{j=1}a_j,\quad a_j=\sum^3_{i=1}c_{j,i}(\frac{T}{T_{ref}})^{i-1}\\
    k_j &= K_{1,j}+K_{2,j}\frac{x_{H2O}}{x_{CO2}}, \quad x_i = \frac{n_i}{n_g}\\
    c_{j,i} &= C_{1,j,i}+C_{2,j,i}\frac{x_{H2O}}{x_{CO2}}+C_{3,j,i}(\frac{x_{H2O}}{x_{CO2}})^2
\end{align}
$T_{ref}$ is 1200 K, and the $K$ and $C$ coefficients are given by the look-up table in \cite{Johanson2011}.
The absorptivity of the gas can be defined as a function of the emissivity \cite{Perry}: 
\begin{align}
    \alpha_g = \epsilon_g(T_s)P_mS_m\sqrt{\frac{T_s}{T_g}},\ S_m = 0.95(2r_c-h)
\end{align}
with $S_m$ being the average path length \cite{Gorog1981}, and $P_m$ being the partial pressure of \ce{CO2} and \ce{H2O}. 

\subsection{Stoichiometry and kinetics}
The reactions occurring in the kiln are described by their reaction rate $r_j = r(T,P,C)$ respectively. The production rates $R$ of each composite are related to the reaction rates:
\begin{align}
    \begin{bmatrix}
        R_s\\R_g
    \end{bmatrix}&= S^Tr
\end{align}
where $R_s$ is the production rate vector of the solids: \ce{CaCO3}, \ce{CaO}, \ce{SiO2}, \ce{AlO2}, \ce{FeO2}, \ce{C_2S}, \ce{C_3S}, \ce{C_3A}, and \ce{C_4AF}, and $R_g$ is the production rate vector of the gasses: $CO_2$, $N_2$ $O_2$, $Ar$, $CO$, $H_2$ $H_2O$, and $C_{sus}$ (suspended carbon). $S$ is the stoichiometry matrix describing the clinker reactions:
\begin{subequations}
\begin{align}
   \text{$r_1$: }& & \ce{CaCO3} &\rightarrow \ce{CaO} + \ce{CO2}\\
   \text{$r_2$: }& & 2\ce{CaO} + \ce{SiO_2} &\rightarrow \ce{C_2S}\\
    \text{$r_3$: }& &\ce{CaO} + \ce{C_2S}&\rightarrow \ce{C_3S}\\
    \text{$r_4$: }& &3\ce{CaO} + \ce{Al_2O_3}&\rightarrow \ce{C_3A}\\
    \text{$r_5$: }& &4\ce{CaO} + \ce{Al_2O_3} + \ce{Fe_2O_3}&\rightarrow \ce{C_4AF}
    \end{align}
\end{subequations}
and the fuel reactions:    
\begin{subequations}
    \begin{align}
   \text{$r_6$: }& & 2\ce{CO}+\ce{O_2}&\rightarrow 2\ce{CO_2}\\
   \text{$r_{7}$: }& & \ce{CO} + \ce{H_2O} &\rightarrow \ce{CO_2} + \ce{H_2}\\
   \text{$r_{8}$: }& & 2\ce{H_2}+\ce{O_2}&\rightarrow 2\ce{H_2O}\\
   \text{$r_9$: }& &2\ce{C} +\ce{ O_2} &\rightarrow 2\ce{CO}\\
   \text{$r_{10}$: }& &\ce{C} + \ce{H_2O}&\rightarrow \ce{CO} + \ce{H_2}\\
   \text{$r_{11}$: }& &\ce{ C} + \ce{CO_2}&\rightarrow 2\ce{CO}.
\end{align}
\end{subequations}
A typical expression for the rate function $r(T,P,C)$ of each reaction is:
\begin{align}
    r_j = k_r(T)\prod_lP_l^{\beta_l}C_l^{\alpha_l},\quad k(T) = k_rT^ne^{-\frac{E_{A}}{RT}}
\end{align}
where $C$ is expressed in Liters, $k(T)$ is the modified Arrhenius equation, $\alpha_l$ is either the stoichiometric or experimental-based values and $\beta_l$ is the power of the partial pressure $P_l = \frac{C_l}{\sum_j C_j}P$. Tables \ref{tab:reaction} and \ref{tab:reaction2} show the reaction coefficients found in the literature for the clinker and fuel reactions respectively.
\begin{table}
    \caption{Reaction rate coefficients of solids in literature.}
    \begin{tabular}{c|c  c | c  c | c c c }%| c}
    & $k_r$ \cite{Spang}& $k_r$ \cite{Mastorakos1999CFDPF} & $E_{A}$ \cite{Spang}  & $E_{A}$ \cite{Mastorakos1999CFDPF} & $\alpha_1$ & $\alpha_2$ & $\alpha_3$ \\%& $\Delta H$\cite{Mujumdar2006}\\
    \hline
    $r_1$ & $1.64\cdot 10^{35}$ & $10^{8}$ & $804.8$& $175.7$& 1& & \\%& -179.4\\
    $r_2$ & $14.8\cdot 10^{8}$&$10^{7}$ & $193.1$& $240$& 2& 1& \\%& 127.6\\
    $r_3$ & $4.8\cdot 10^{8}$& $10^{9}$& $255.9$& $420$& 1& 1& \\%& -16\\
    $r_4$ & $300\cdot 10^{8}$&$10^{8}$ & $193.8$& $310$& 3& 1& \\%& -21.8\\
    $r_5$ & $30\cdot10^{11}$& $10^{8}$& $184.9$ & $330$& 4& 1&1 \\%& 41.3\\
    \hline
    \end{tabular}    
    \footnotesize{The units of the reactions are $[\frac{1}{\text{hr}}]$ and $[\frac{\text{kg}}{\text{m}^3\text{s}}]$ respectively for \cite{Spang} and \cite{Mastorakos1999CFDPF}, with the activation energy $E_{A}$ given in $[\frac{\text{kj}}{\text{mol}}]$. All $\beta$ and $n$ are zero.}
       \label{tab:reaction}
\end{table}

\begin{table}[tb]
    \caption{Reaction rate coefficients of gasses in literature}\label{tab:reaction2}%
    \begin{tabular*}{0.5\textwidth}{c | c| c|c | c | c |c | c }%| c   }
    & Unit & $k_r $ & n  & $E_{A}$ & $\alpha_1$ & $\alpha_2$ & $\beta_2$ \\%& $\Delta H$\cite{BASU2018211}\\
    \hline
        $r_6$ \cite{Guo2003}&$\frac{\text{kg}}{\text{m}^3 \cdot\text{s}}$&$7.0\cdot10^4$ & 0 & $66.51$& $1$\footnotemark[1]  & $1$\footnotemark[1]  &0 \\%& 284\\
        $r_7$ \cite{JONES1988}&$\frac{\text{mol}}{\text{m}^3 \cdot\text{s}}$& $2.75\cdot10^6$ & 0 & 83.68 & 1& 1& 0 \\%& 41.2\\
        $r_8$ \cite{Karkach1999}&$\frac{\text{mol}}{\text{m}^3 \cdot\text{s}}$& $1.37\cdot10^{6}$ & 0.51 & $295.48$ & 1 & 1 & 0 \\%& 242\\
        $r_9$  \cite{Walker1985}& $\frac{\text{mol}}{\text{m}^3 \cdot\text{s}}$& $8.82\cdot10^{13}$ & 0 &239 & $1$\footnotemark[1]  & $0.5$\footnotemark[1]  & 0 \\%& 111\\        
        $r_{10}$ \cite{BASU2018211}& $\frac{\text{mol}}{\text{m}^3 \cdot\text{s}}$& $2.62\cdot10^8$ & 0 & 237 & 1 & 0 &  0.57 \\%& -131\\
        $r_{11}$ \cite{BASU2018211}& $\frac{mol}{m^3 s}$& $3.1\cdot10^6$ & 0 & 215 & 1 & 0  & 0.38 \\%& -172\\
     \hline
    \end{tabular*}
    \footnotesize{All $\beta_1$ is zero, the unit of the activation energy $E_{A}$ is $[\frac{kJ}{mol}]$. $^1$ is unclear in source.}
\end{table}

\subsection{Mass balance}
Combining the material flux and reaction rates, the mass balances are given in their concentration form:
%The mass balances are given in their concentration forms. The change in concentrations depends on the material flux and reaction rates:
\begin{align}
    \partial_t\hat{C}_{i,s} &= -\partial_zN_{i,s} + R_{s,i}\\
    \partial_t\hat{C}_{i,g} &= -\partial_zN_{i,g} + R_{g,i}.
\end{align}

\subsection{Energy balance}
The energy balances are obtained from the heat transport and the enthalpy flux $\Tilde{H}_i=H(T_i,P,N_i)$:
%The energy balances are given by their energy densities. The change in energy density is formulated through the enthalpy flux ($\Tilde{H}_i=H(T_i,P,N_i)$) and heat transport:
{\small
\begin{align}
    \partial_t\hat{U}_s &= -\partial_z(\Tilde{H}_s + \Tilde{Q}_s) + \frac{Q^{rad}_{gs}+Q^{rad}_{ws} + Q^{cv}_{gs} + Q^{cv}_{ws}}{V_{\Delta}} - J_{sg}\\
    %\partial_t\hat{U}_s &= -\partial_z\Tilde{H}_s -\partial_z\hat{V}_s\Tilde{Q}_s + \frac{Q^{rad}_{gs}+Q^{rad}_{ws} + Q^{cv}_{gs} + Q^{cv}_{ws}}{V_{\Delta}}\\       
    \partial_t\hat{U}_g &= -\partial_z(\Tilde{H}_g+\Tilde{Q}_g) - \frac{Q^{rad}_{gs}+Q^{rad}_{gw}+Q^{cv}_{gs}+Q^{cv}_{gw}}{V_{\Delta}} + J_{sg}\\
    \partial_t\hat{U}_w &= -\partial_z\Tilde{Q}_w+ \frac{Q^{rad}_{gw}- Q^{rad}_{ws} + Q^{cv}_{gw} - Q^{cv}_{ws}}{V_{\Delta w}} .
\end{align}
}%
where $J_{sg}=H(T_s,P,r_1)$ is the phase transition term of reaction $r_1$ for the produced \ce{CO2}.
% {\color{red}
% \subsubsection{Boundary conditions}
% The boundaries consist of the pressures at each end of the kiln, and the temperatures and masses entering the kiln. For the solids, the compound is assumed to have the same temperature and velocity, while the entering gas is a mixed of injected fuel and air.
% Assuming mass conservation through the kiln end and the fuel velocity being controlled, the air flow is given by
% \begin{align}
%     \dot{m}_{g,in} &=\rho_gA_gv_g = \rho_{air}A_{air}v_{air} + \rho_fA_fv_f\\
%     A_g &= A_{air}+A_f
% \end{align}
% with the gas temperature entering the kiln given by
% \begin{align}
%     \dot{m}_{g,in}T_{g,in} = \dot{m}_{f}T_{g,f} + \dot{m}_{air}T_{air}
% \end{align}
%  a weighted sum of air and fuel temperatures, weighted by the mass flows.
% }
\subsection{Algebraic equations}
The algebraic part consists of the volume density relation:
\begin{align}
    \hat{V}_{g} + \hat{V}_{s} = \hat{V}_{\Delta} = 1
    %\sum_i \Delta P_i = P_0-P_N
\end{align}
to constrain the pressures, the energy density relations in \eqref{eq:EnergyAlgebra} for the temperatures, and geometry of \eqref{eq:geoAlg} for the fill angles.

% In addition to the differential descriptions, the model includes the algebraic relations for the volume density:
% \begin{align}
%     \hat{V}_{g} + \hat{V}_{s} = \hat{V}_{\Delta} = 1
%     %\sum_i \Delta P_i = P_0-P_N
% \end{align}
% and the energy density relations given by (\ref{eq:EnergyAlgebra}).
% $\theta$ based on the geometry of \eqref{eq:geoAlg}
\subsubsection{Boundary conditions}
The boundary of the model is given as the pressure at the beginning of the kiln (direction of gas flow), and the temperatures and influx of solids and gasses at their respective kiln ends.
\section{Data and Evaluation}
%The mathematical model discussed has a reliance on the availability of data on the properties of the different compounds and materials present in the kiln. These properties has to either be experimentally evaluated or based on available literature. In the literature, the properties are as given in table \ref{tab:Data-Coeff-solid} and table \ref{tab:Data-Coeff-gas} for the solid and gas respectively. Table \ref{tab:heatCap} includes the molar heat capacity of the compounds in the kiln: $c_m = C_0 + C_1T + C_2T^2$.
The model relies on physical properties, provided either from experiments or available through the literature. 
Based on the literature, table \ref{tab:Data-Coeff-solid} and table \ref{tab:Data-Coeff-gas} provide the material properties for the solids and gasses respectively. 
Table \ref{tab:heatCap} includes the molar heat capacity $c_m = C_0 + C_1T + C_2T^2$.
\begin{table}
    \centering
    \caption{Material properties of the solid phase}% - perry might be profiles}
    \begin{tabular}{c|c|c|c|c}
    \hline
           &\shortstack{Thermal\\ Conductivity} & Density & \shortstack{Molar \\mass} & \shortstack{Enthalpy of \\formation}\\ \hline
         Units    & $\frac{\text{W}}{\text{K}\cdot\text{m}}$ & $\frac{\text{g}}{\text{cm}^3}$ & $\frac{\text{g}}{\text{mol}}$ & $\frac{\text{kJ}}{\text{mol}}$ \\ \hline
         $CaCO_3$ &  2.248$^a$& 2.71$^b$  &100.09$^b$ & - 1207.6$^b$\\ 
         $CaO$     & 30.1$^c$ &  3.34$^b$ &56.08$^b$ & -634.9$^b$\\ 
         $SiO_2$   &  1.4$^{a,c}$& 2.65$^b$  &60.09$^b$ & -910.7$^b$\\ 
         $Al_2O_3$ &  12-38.5$^c$ 36$^a$& 3.99$^b$  &101.96$^b$& -1675.7$^b$\\ 
         $Fe_2O_3$ &  0.3-0.37$^c$& 5.25$^b$  &159.69$^b$&	-824.2$^b$\\ \hline
         $C2S$     &  3.45$\pm$0.2$^d$& 3.31$^d$  &$172.24^g$& -2053.1$^h$\\ 
         $C3S$     & 3.35$\pm$0.3$^d$ & 3.13$^d$ & 228.32$^b$&-2704.1$^h$\\ 
         $C3A$     &  3.74$\pm$0.2$^e$& 3.04$^b$ & 270.19$^b$&-3602.5$^h$\\ 
         $C4AF$    &  3.17$\pm$0.2$^e$& 3.7-3.9$^f$  &$485.97^g$ & -4998.6$^h$\\ \hline        
    \end{tabular}   
    
    \footnotesize{$^a$ Data from \cite{Perry}, $^b$ Data from \cite{CRC2022}, $^c$ Data from \cite{Ichim2018}, $^d$ Data from \cite{PhysRevApplied},\\ $^e$ Data from \cite{Du2021}, $^f$ Data from \cite{Portland}, $^g$ Computed from the above results,\\ $^h$ Computed from \cite{Mujumdar2006} and Hess' Law.} 
    \label{tab:Data-Coeff-solid}
\end{table}
\begin{table}
    \centering
    \caption{Material properties of the gas phase}
    \begin{tabular}{c|c|c|c|c|c}
    \hline
           &\shortstack{Thermal\\ Conduc-\\tivity$^a$} & \shortstack{Molar\\ mass$^a$} & Viscosity$^a$ & \shortstack{diffusion\\ Volume$^b$}& \shortstack{Enthalpy\\ of \\formation$^f$}\\ \hline
        Units   & $\frac{10^{-3}\text{W}}{\text{K}\cdot\text{m}}$ & $\frac{\text{g}}{\text{mol}}$ & $\mu \text{Pa}\cdot\text{s}$ & \text{cm}$^3$ & $\frac{\text{kJ}}{\text{mol}}$\\ \hline
         \ce{CO_2}    &\shortstack{\strut 16.77$^c$\\ 70.78$^e$ }  & 44.01  & \shortstack{\strut 15.0$^c$\\ 41.18$^e$ }& 16.3 & - 395.5\\\hline
         \ce{N_2} & \shortstack{\strut 25.97$^c$\\  65.36$^e$ }  &28.014 &  \shortstack{ \strut 17.89$^c$\\  41.54$^e$ } & 18.5 & 0\\\hline
         
         \ce{O_2}  & \shortstack{\strut 26.49$^c$\\  71.55$^e$ } &  31.998 &  \shortstack{\strut 20.65$^c$\\ 49.12$^e$ }&  16.3 & 0\\\hline
         
         \ce{Ar}  & \shortstack{\strut 17.84$^c$\\ 43.58$^e$ }& 39.948&  \shortstack{\strut  22.74$^c$\\  55.69$^e$ }&  16.2 & 0\\\hline
         
         \ce{CO}  & \shortstack{\strut 25$^c$\\  43.2$^d$ } & 28.010& \shortstack{\strut 17.8$^c$\\  29.1$^e$ } &  18& -110.5\\\hline
         
         \ce{C_{sus}} & - & 12.011 & - &  15.9 & 0\\\hline
         
         \ce{H_2O} & \shortstack{\strut 609.50$^c$\\  95.877$^e$ }&18.015 &  \shortstack{\strut 853.74$^c$\\  37.615$^e$ }&   13.1 & -241.8\\\hline
         
         \ce{H_2} & \shortstack{\strut 193.1$^c$\\ 459.7$^e$ }&2.016 &  \shortstack{ \strut 8.938$^c$\\ 20.73$^e$ }& 6.12 & 0\\\hline
    \end{tabular}

    $^a$ Data from \cite{CRC2022}, $^b$ Data from \cite{Poling2001Book}, $^c$ Value at 300 K, $^d$ Value at 600 K,\\ $^e$ Value at 1000 K, $^f$ Data from \cite{BASU2018211}.
    \label{tab:Data-Coeff-gas}
\end{table}
In a real kiln, several more reactions occur with compounds not included here, e.g. hydrocarbon fuel or alternative fuels.
The model's flexibility allows for easy extensions of compounds, by
adding mass balances and the needed reactions.
The layout also allows for easy replacing a specific model part, e.g. a reaction rate, with a desired formulation.

\begin{table}
    \centering
    \caption{Molar heat capacity}
    \begin{tabular}{c|c |c | c| c}
            & $C_0$ & $C_1$ & $C_2$ & Temperature range\\ \hline
         Units & $\frac{\text{J}}{\text{mol}\cdot\text{K}}$& $\frac{10^{-3}\text{J}}{\text{mol}\cdot\text{K}^2}$&$\frac{10^{-5}\text{J}}{\text{mol}\cdot\text{K}^3}$& K\\\hline
         \ce{CaCO_3}$^a$ &  23.12 &  263.4 & -19.86& 300 - 600\\ 
         \ce{CaO}$^b$     &  71.69& -3.08  & 0.22  & 200 - 1800\\ 
         \ce{SiO_2}$^b$   & 58.91 &  5.02 & 0& 844 - 1800\\ 
         \ce{Al_2O_3}$^b$ &  233.004&-19.5913   &0.94441  & 200 - 1800\\ 
         \ce{Fe_2O_3}$^a$ & 103.9  & 0 & 0 & -\\ \hline
         \ce{C2S}$^b$     &  199.6& 0  &0  & 1650 - 1800\\ 
         \ce{C3S}$^b$     &  333.92&  -2.33&  0& 200 - 1800\\ 
         \ce{C3A}$^b$    & 260.58  & 9.58/2 & 0  &298 - 1800\\ 
         \ce{C4AF}$^b$  &  374.43& 36.4 & 0 &  298 - 1863\\ \hline
         \ce{CO_2}$^a$    & 25.98 &43.61 &-1.494 & 298 - 1500\\
         \ce{N_2}$^a$ & 27.31&5.19 &-1.553e-04 & 298 - 1500\\ 
         \ce{O_2}$^a$ & 25.82&12.63 &-0.3573  & 298 - 1100\\ 
         \ce{Ar}$^a$ & 20.79 & 0 & 0  & 298 - 1500\\ 
         \ce{CO}$^a$ & 26.87& 6.939  & -0.08237  & 298 - 1500\\ 
         \ce{C_{sus}}$^a$& -0.4493& 35.53 & -1.308 & 298 - 1500\\ 
         \ce{H_2O}$^a$&30.89 & 7.858  &0.2494 & 298 - 1300\\ 
         \ce{H_2}$^a$&   28.95& -0.5839&  0.1888 & 298 - 1500\\ \hline         
    \end{tabular}    
    \footnotesize{ $^a$ Based on data from  \cite{CRC2022}, $^b$ Coefficients from \cite{HANEIN2020106043}. }
    \label{tab:heatCap}
\end{table}

\section{Simulation}
To demonstrate the model, we consider a 50 m long kiln with a 2 m inner radius and 2\% inclination, split into 10 segments.
Matlab's ode15s is used to simulate the 200 differential equations and 40 algebraic equations.
Initially, the kiln is filled with air with 1\% \ce{H2O}, the air and wall temperatures are 800$^\circ$C along the kiln, and a pressure profile of 1.00005 bar to 1.00010 bar induces a slight flow through the kiln.

To operate the kiln at a 13\% fill ratio (common load) at 4 rpm, the solid inflow is 102 ton/h at 800$^\circ$C with a composition of 73\% \ce{CaO}, 22.5\% \ce{SiO2}, 3\% \ce{Al2O2}, and 1.5\% \ce{Fe2O3}, giving an influx of 778 mol/m$^3$ at 0.048 m/s.
For the heating a 300 kcal/kg clinker is used, given a fuel inflow of 3.9 ton/h carbon at 1200$^\circ$C, an influx of 2.4 mol/m$^3$ at 3 m/s.
The secondary air influx is 7.2 mol/m$^3$ at 6 m/s and 1200$^\circ$C.
For the reaction rates we use Mastorakas' values \cite{Mastorakos1999CFDPF}, and for efficient computation, the sgn(x) in \eqref{eq:vel} is approximated using tanh(30x).

\subsection{Tuning}
To achieve a behavior resembling the real world, the following parts were tuned visually. The solid densities by a factor of $1/9$; to match the nominal fill ratio and nominal throughput (ton/h) for the dimension used. 
Based on the exit solid/gas composition, the reactions $r_1$-$r_6$ were tuned by the factors 5, 200, 60, $5\times 10^3$, $5\times10^6$, and $3\times10^4$ respectively.

\subsection{Performance}
\begin{figure}
    \centering
    \includegraphics[width=0.45\textwidth,trim={1.20cm, 0.05cm, 1.75cm, 0.35cm},clip]{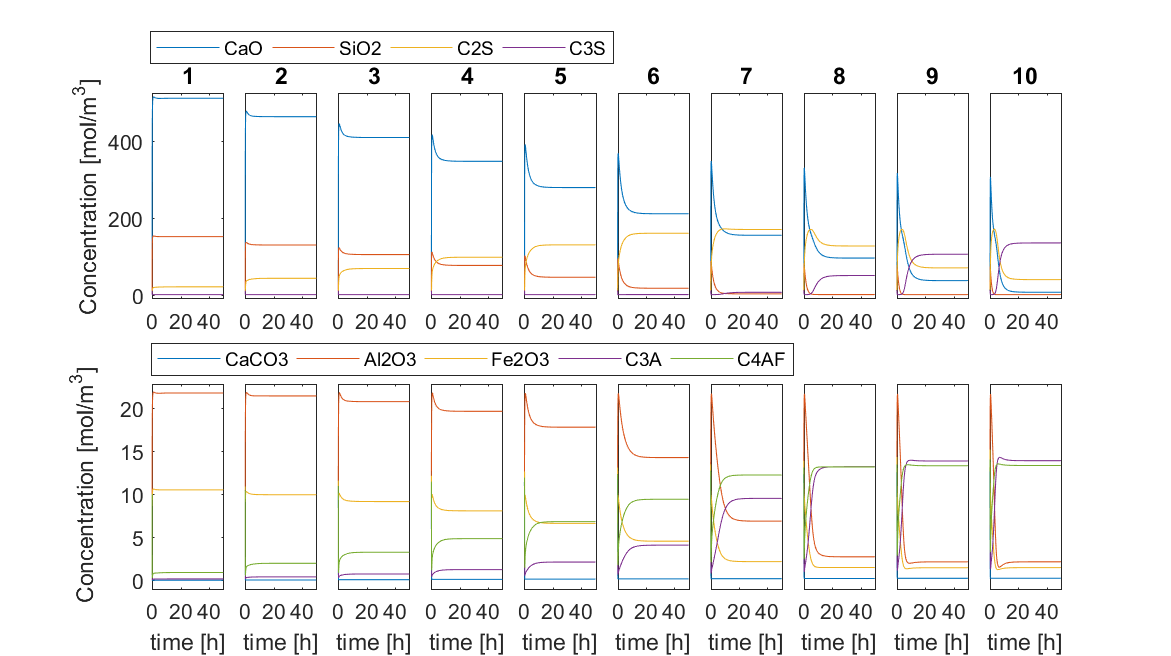}  
    \caption{Dynamic evolution of solid concentrations along the kiln for a 50-hour simulation, showing the transition and steady state of nine compounds in the clinker.}
    \label{fig:solids}
\end{figure}

\begin{figure}
    \centering
    \includegraphics[width=0.45\textwidth,trim={1.1cm, 0.05cm, 1.75cm, 0.35cm},clip]{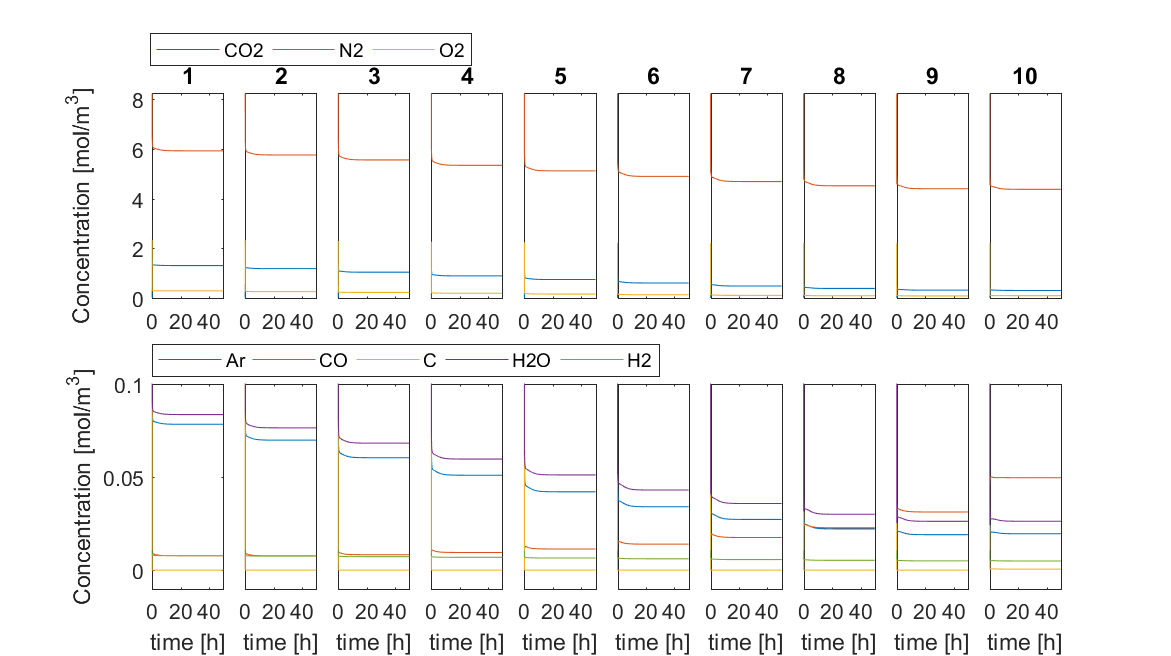}  
    \caption{Dynamic evolution of gas concentrations along the kiln for a 50-hour simulation, showing the transition and steady state of 8 compounds in the gas.}
    \label{fig:gasses}
\end{figure}

\begin{figure}
    \centering
    \includegraphics[width=0.45\textwidth,trim={1.16cm, 0.05cm, 1.85cm, 0.95cm},clip]{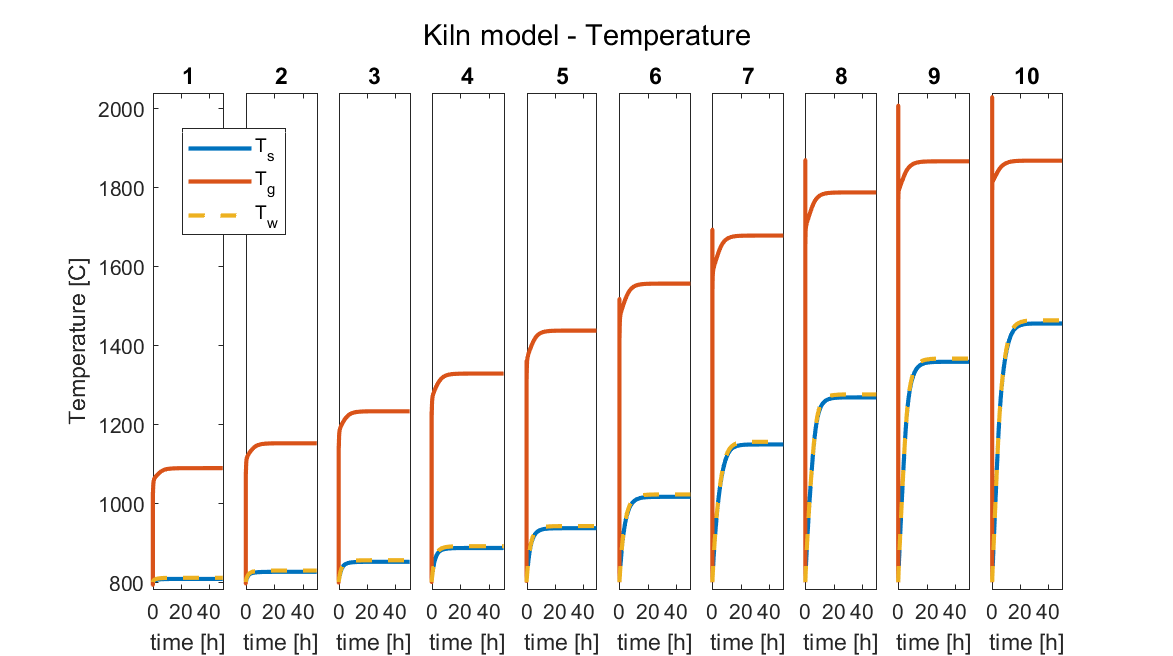}  
    \caption{Dynamic evolution of temperature along the kiln for a 50-hour simulation, showing the bulk temperatures of the solid, gas, and wall phases. }
    \label{fig:temps}
\end{figure}

\begin{figure}
    \centering
    \includegraphics[width=0.5\textwidth,trim={0.35cm, 0.05cm, 0.15cm, 0.25cm},clip]{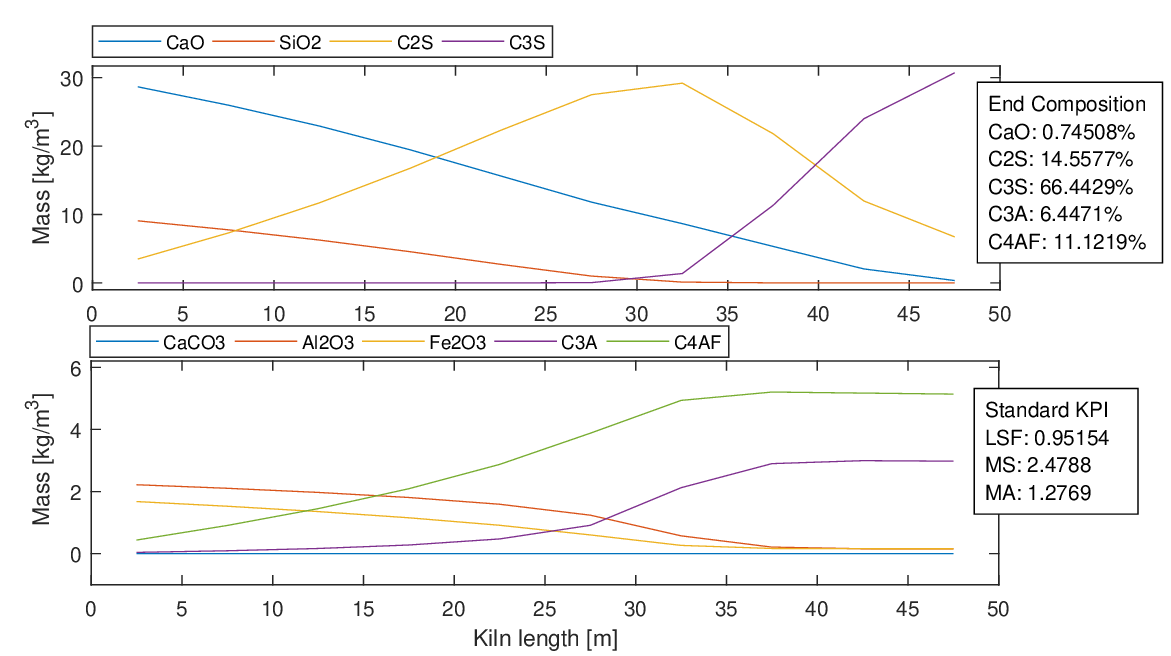}    
    \caption{Steady-state mass concentrations of solid compounds along the kiln.
Typical quality KPIs for the raw meal entering the kiln, as well as the resulting compound percentages in the mixture at the end of the kiln, are included.}
    \label{fig:solidstatic}
\end{figure}

For a 50-hour simulation, fig. \ref{fig:solids}-\ref{fig:temps} shows the evolution in compound concentrations and temperatures. 
We can see the model settles in around 28 hours, resembling the usual rotary kiln start-up time of 24-48 hours. From the solid compounds, we see how \ce{C2S} gets formed from \ce{CaO} becoming the dominating compound until \ce{C3S} starts forming from the \ce{C2S} at around 1300 $^\circ$C. In the gas, we see the fuel $C_s$ is consumed immediately giving an increase in \ce{CO2}, a decrease in \ce{O2}, and a sudden increase followed by a decline in \ce{CO}, resulting in 17.0\%, 3.8\%, and 0.1\% of the gas outlet respectively.
The observed temperature profiles likewise resemble that of a typical rotary kiln, with solid temperature starting around 800$^\circ$C and reaching above 1450$^\circ$C, while the gas temperature ranges from 1870$^\circ$C at the end to 1090$^\circ$C at the beginning of the kiln. A steady-state pressure drop of $~27$ Pa along the kiln, giving an exit velocity of 5.6 m/s.
Fig. \ref{fig:solidstatic} shows the steady-steady profile of solids in mass concentration. The standard Key Performance Indicator (KPI) measures of clinker quality are given (LSF, MS, MA). These are based on the input raw meal, while the solid composition at the kiln end is given in percentages: 66.4\% \ce{C3S}, 14.6\% \ce{C2S}, 6.4\% \ce{C3A}, 11.1\% \ce{C4AF}, and 0.7\% \ce{CaO}. The seen percentages are within the typical  range of cement clinker, with the given KPI.

In general, we observe the model follows the expected behaviors for a cement rotary kiln, with only slight tuning needed; e.g. reaction activation energy to make reaction starting temperatures more accurate.

\section{Conclusion}
In this paper, we presented a dynamic model of a rotary kiln for cement clinker production based on first engineering principles. The modeling focuses on the formation of clinker compounds. The model is described systematically, outlining the used thermophysical, chemical, and transportation models; enabling an easy overview of the dynamic details included, and how to extend it. 
The paper includes a collection of necessary property data of the different materials from the literature. The model was illustrated with a simulation, showcasing the dynamics and steady states of the kiln, resembling the expected behavior.
For existing cement plants, the proposed simulation model can be used as a simplifying step in the development of advance process control systems.

%The paper presented the mathematical description of the model, together with 500 hour simulation scenario and a collection of necessary property data of the different materials, found in the literature.

%\subsection{where to go}

%The pre-heating tower and calciner were introduced in modern production setups, to shorten the kiln and speed up the reactions.
%further work could model these and the cooler, to combine the entire dynamic profile of the clinker production.

% while the model is detailed descriped at some points, assumptions are used at others, the model formulation is deliberately kept expanded so the details of each part is clear, but also so that the model can be updated by easily by exchanging only the specific part and keeping the rest intact. which allows for easy simplification or expansion of complexity.

% improvements: separate fuel flow dynamics?, expanded reaction/composite list
% 3-D description
% loss to environment.

% devolation, flame distribution/placement - same air, cfd
% radiation: axial components
% particle aspect
% kiln coating
% turbulence
% false air
% the other sections in the process line
%%%%%%%%%%%%%%%%%%%%%%%%%%%%%%%%%%%%%%%%%%%%%%%%%%%%%%%%%%%%%%%%%%%%%%%%%%%%%%%%

%\section*{APPENDIX}

\section*{ACKNOWLEDGMENT}
The authors would like to acknowledge the financial support by the Innovation Fund of Denmark (No. 053-00012B – Industrial Postdoc Project: Green-Digital Transition in Cement Production)

\bibliographystyle{IEEEtran}
\bibliography{IEEEabrv,biblio}

\begin{thebibliography}{10}
\providecommand{\url}[1]{#1}
\csname url@rmstyle\endcsname
\providecommand{\newblock}{\relax}
\providecommand{\bibinfo}[2]{#2}
\providecommand\BIBentrySTDinterwordspacing{\spaceskip=0pt\relax}
\providecommand\BIBentryALTinterwordstretchfactor{4}
\providecommand\BIBentryALTinterwordspacing{\spaceskip=\fontdimen2\font plus
\BIBentryALTinterwordstretchfactor\fontdimen3\font minus
  \fontdimen4\font\relax}
\providecommand\BIBforeignlanguage[2]{{%
\expandafter\ifx\csname l@#1\endcsname\relax
\typeout{** WARNING: IEEEtran.bst: No hyphenation pattern has been}%
\typeout{** loaded for the language `#1'. Using the pattern for}%
\typeout{** the default language instead.}%
\else
\language=\csname l@#1\endcsname
\fi
#2}}

\bibitem{CO2Techreport}
J.~Lehne and F.~Preston, ``Making concrete change: Innovation in low-carbon
  cement and concrete,'' Chatham House, Tech. Rep., june 2018.

\bibitem{Mujumdar2006}
K.~S. Mujumdar and V.~V. Ranade, ``Simulation of rotary cement kilns using a
  one-dimensional model,'' \emph{Chemical Engineering Research and Design},
  vol.~84, pp. 165--177, 2006.

\bibitem{Hanein2017}
T.~Hanein, F.~P. Glasser, and M.~N. Bannerman, ``One-dimensional steady-state
  thermal model for rotary kilns used in the manufacture of cement,''
  \emph{Adv. Appl. Ceram.}, vol. 116, no.~4, pp. 207--215, 2017.

\bibitem{Mastorakos1999CFDPF}
E.~Mastorakos, A.~Massias, C.~D. Tsakiroglou, D.~A. Goussis, V.~N. Burganos,
  and A.~C. Payatakes, ``Cfd predictions for cement kilns including flame
  modelling, heat transfer and clinker chemistry,'' \emph{Applied Mathematical
  Modelling}, vol.~23, pp. 55--76, 1999.

\bibitem{Spang}
H.~A. Spang, ``A dynamic model of a cement kiln,'' \emph{Automatic}, vol.~8,
  pp. 309--323, 1972.

\bibitem{SUN2020}
C.~Sun, J.~Zhao, S.~Li, and P.~Jiang, ``First-principle modeling and simulation
  of cement rotary kiln,'' in \emph{2020 Chinese Control And Decision
  Conference (CCDC)}, 2020, pp. 3267--3272.

\bibitem{LIU2016}
H.~Liu, H.~Yin, M.~Zhang, M.~Xie, and X.~Xi, ``Numerical simulation of particle
  motion and heat transfer in a rotary kiln,'' \emph{Powder Technology}, vol.
  287, pp. 239--247, 2016.

\bibitem{GINSBERG2011}
T.~Ginsberg and M.~Modigell, ``Dynamic modelling of a rotary kiln for
  calcination of titanium dioxide white pigment,'' \emph{Computers \& Chemical
  Engineering}, vol.~35, no.~11, pp. 2437--2446, 2011.

\bibitem{CircSeg}
\BIBentryALTinterwordspacing
E.~W. Weisstein. (2023, feb.) Circular segment. MathWorld--A Wolfram Web
  Resource. [Online]. Available:
  \url{https://mathworld.wolfram.com/CircularSegment.html}
\BIBentrySTDinterwordspacing

\bibitem{Perry}
D.~W. Green and R.~H. Perry, Eds., \emph{Perry's Chemical Engineers' Handbook},
  8th~ed.\hskip 1em plus 0.5em minus 0.4em\relax McGraw Hill, 2008.

\bibitem{Saeman1951}
W.~C. Saeman, ``Passage of solids through rotary kilns,'' \emph{Chemical
  Engineering Progress}, vol.~47, no.~10, pp. 508--514, 1951.

\bibitem{Yang2003}
R.~Yang, R.~Zou, and A.~Yu, ``Microdynamic analysis of particle flow in a
  horizontal rotating drum,'' \emph{Powder Technology}, vol.~30, 2003.

\bibitem{Yamane1998}
K.~Yamane, M.~Nakagawa, S.~Altobelli, T.~Tanaka, and Y.~Tsuji,
  ``\BIBforeignlanguage{English}{Steady particulate flows in a horizontal
  rotating cylinder},'' \emph{\BIBforeignlanguage{English}{Physics of Fluids}},
  vol.~10, no.~6, pp. 1419--1427, June 1998.

\bibitem{MortenPHD}
M.~N. Pedersen, ``Co-firing of alternative fuels in cement kiln burners,''
  Ph.D. dissertation, Technical University of Denmark, 2018.

\bibitem{Darcy-Howel}
G.~W. Howell and T.~M. Weathers, \emph{Aerospace Fluid Component Designers'
  Handbook. Volume I, Revision D}.\hskip 1em plus 0.5em minus 0.4em\relax TRW
  Systems Group, 1970.

\bibitem{HESSELGREAVES20171}
J.~E. Hesselgreaves, R.~Law, and D.~A. Reay, ``Chapter 1 - introduction,'' in
  \emph{Compact Heat Exchangers}, 2nd~ed.\hskip 1em plus 0.5em minus
  0.4em\relax Butterworth-Heinemann, 2017.

\bibitem{Poling2001Book}
B.~E. Poling, J.~M. Prausnitz, and J.~P. O'Connel, \emph{The Properties of
  Gases and Liquids}.\hskip 1em plus 0.5em minus 0.4em\relax McGraw-Hill, 2001.

\bibitem{Tscheng1979}
S.~H. Tscheng and A.~P. Watkinson, ``Convective heat transfer in a rotary
  kiln,'' \emph{Can. J. Chem. Eng.}, vol.~57, pp. 433--443, 1979.

\bibitem{Sutherland1893}
W.~Sutherland, ``Lii. the viscosity of gases and molecular force,''
  \emph{Philosphical Magazine series 5}, vol.~36, no. 223, pp. 507--531, 1893.

\bibitem{Wilke1950}
C.~R. Wilke, ``A viscosity equation for gas mixtures,'' \emph{The Journal of
  Chemical Physics}, vol.~18, no.~4, pp. 517--519, 1950.

\bibitem{Johanson2011}
R.~Johansson, B.~Leckner, K.~Andersson, and F.~Johnsson, ``Account for
  variations in the \ce{H2O} to \ce{CO2} molar ratio when modelling gaseous
  radiative heat transfer with the weighted-sum-of-grey-gases model,''
  \emph{Combustion and Flame}, vol. 158, pp. 893--901, 2011.

\bibitem{Gorog1981}
J.~P. Gorog, J.~K. Brimacombe, and T.~N. Adams, ``Radiative heat transfer in
  rotary kilns,'' \emph{Metall Mater Trans B}, vol.~12, 1981.

\bibitem{Guo2003}
Y.~Guo, C.~Chan, and K.~Lau, ``Numerical studies of pulverized coal combustion
  in a tubular coal combustor with slanted oxygen jet,'' \emph{Fuel}, vol.~82,
  pp. 893--907, 2003.

\bibitem{JONES1988}
W.~P. Jones and R.~P. Lindstedt, ``Global reaction schemes for hydrocarbon
  combustion,'' \emph{Combustion and Flame}, vol.~73, 1988.

\bibitem{Karkach1999}
S.~P. Karkach and V.~I. Osherov, ``Ab initio analysis of the transition states
  on the lowest triplet \ce{H2O2} potential surface,'' \emph{The Journal of
  Chemical Physics}, vol. 110, no.~24, pp. 11\,918--11\,927, 1999.

\bibitem{Walker1985}
P.~L. Walker, \emph{Char Properties and Gasification}.\hskip 1em plus 0.5em
  minus 0.4em\relax Springer Netherlands, 1985, pp. 485--509.

\bibitem{BASU2018211}
P.~Basu, ``Chapter 7 - gasification theory,'' in \emph{Biomass Gasification,
  Pyrolysis and Torrefaction}, 3rd~ed.\hskip 1em plus 0.5em minus 0.4em\relax
  Academic Press, 2018.

\bibitem{CRC2022}
J.~Rumble, Ed., \emph{CRC handbook of chemistry and physics}, 103rd~ed.\hskip
  1em plus 0.5em minus 0.4em\relax CRC Press, 2022.

\bibitem{Ichim2018}
A.~Ichim, C.~Teodoriu, and G.~Falcone, ``Estimation of cement thermal
  properties through the three-phase model with application to geothermal
  wells,'' \emph{Energies}, vol.~11, no.~10, 2018.

\bibitem{PhysRevApplied}
A.~Qomi, M.~Javad, F.-J. Ulm, and R.~J.-M. Pellenq, ``Physical origins of
  thermal properties of cement paste,'' \emph{Phys. Rev. Appl.}, vol.~3, Jun
  2015.

\bibitem{Du2021}
Y.~Du and Y.~Ge, ``Multiphase model for predicting the thermal conductivity of
  cement paste and its applications,'' \emph{Materials}, vol.~14, no.~16, 2021.

\bibitem{Portland}
G.~C. Bye, Ed., \emph{Portland Cement: Composition, Production and Properties},
  2nd~ed.\hskip 1em plus 0.5em minus 0.4em\relax Thomas Telford, 1999.

\bibitem{HANEIN2020106043}
T.~Hanein, F.~P. Glasser, and M.~N. Bannerman, ``Thermodynamic data for cement
  clinkering,'' \emph{Cem Concr Res}, vol. 132, 2020.

\end{thebibliography}

\end{document}